\documentclass[11pt]{article}
\usepackage{times, amsfonts, amsmath, amssymb, latexsym, setspace, epsfig}
\usepackage{subfigure, graphics}
\graphicspath{{figs/}}
\usepackage{multirow}
\usepackage{algorithm,algorithmic}
\usepackage{array}
\usepackage{authblk}
\newcolumntype{P}[1]{>{\centering\arraybackslash}p{#1}}

\def\Ex{{\rm I\!E}}
\def\Pr{{\rm I\!P}}
\def\be{\begin{equation}}
\def\ee{\end{equation}}
\def\bea{\begin{eqnarray*}}
\def\eea{\end{eqnarray*}}
\def\bean{\begin{eqnarray}}
\def\eean{\end{eqnarray}}
\def\nn{\nonumber}
\def\nin{\noindent}

\def\ra{\rightarrow}

\def\Bl{\Bigl}
\def\Br{\Bigr}
\def\wt{\widetilde}

\def\R{{\bf R}}

\def\alp{\alpha}

\def\del{\delta}
\def\eps{\epsilon}

\def\th{\theta}
\def\Th{\Theta}

\def\logLRm{{\rm logLR}_m}
\def\logLRmn{{\rm logLR}_{m,n}(\hat{\th}_n)}
\def\logLRI{{\rm logLR}_I}

\newtheorem{Theorem}{Theorem}
\newtheorem{Proposition}{Proposition}
\newtheorem{Corollary}{Corollary}
\newtheorem{Lemma}{Lemma}

\begin{document}
\setstretch{1.1}

\title{Tail bounds for empirically standardized sums}

\author{Guenther Walther\thanks{Research supported by NSF grants DMS-1501767 and DMS-1916074}\\
Department of Statistics, 390 Jane Stanford Way \\
        Stanford University, Stanford, CA 94305 \\
        gwalther@stanford.edu}

\date{March 2022}

\maketitle

\begin{abstract}
Exponential tail bounds for sums play an important role in statistics, but the example of the $t$-statistic
shows that the exponential tail decay may be lost when population parameters need to be estimated from the data.
However, it turns out that if Studentizing is accompanied by estimating
the location parameter in a suitable way, then the $t$-statistic regains the
exponential tail behavior. Motivated by this example, the paper analyzes other ways of empirically
standardizing sums and establishes tail bounds that are sub-Gaussian or even closer to normal for the
following settings:
Standardization with Studentized contrasts for  normal observations,
standardization with the log likelihood ratio statistic for observations from an exponential family,
and standardization via self-normalization for observations from a symmetric distribution
with unknown center of symmetry. 
The latter standardization gives rise to a novel scan statistic for heteroscedastic data whose
asymptotic power is analyzed in the case where the observations have a log-concave distribution.
\end{abstract}

\vfill



\section{Introduction}

Tail bounds and concentration inequalities for sums of independent random variables play a key role in 
statistics and machine learning, see e.g. van der Vaart and Wellner~(1996), Boucheron et al.~(2013), Vershynin~(2018),
or Wainwright~(2019). Of particular importance are exponential tails bounds, which typically involve
the expected value of the sum as well as a scale factor such as the variance. On the other hand, few results
seem to be available when these parameters need to be estimated from the data, as may be required to make
statistical methodology operational.
The most prominent example is the $t$-statistic: If $X_1,\ldots,X_m$ are i.i.d. N($\mu,\sigma^2$), then
\be \label{Tstat}
T\ :=\ \frac{\frac{1}{\sqrt{m}} \sum_{i=1}^m (X_i -\mu)}{
\sqrt{\frac{1}{m-1} \sum_{i=1}^m \bigl(X_i -\overline{X}\bigr)^2}}
\ee
has the heavy algebraic tails of the $t_{m-1}$-distribution, so estimating $\sigma^2$
with the sample variance comes at the expense of losing the exponential tail decay. This paper explores
the case where the expectation $\mu$ is also unknown and must be estimated. This is the typical setting
for scan statistics, where observations in a scan window are assessed against an unknown baseline which is
estimated with the sample mean of all observations, see e.g. Yao~(1993). Corollary~\ref{tstatistic2} below
shows that, rather than
exacerbating the situation, this additional estimation step actually restores the sub-Gaussian tail bound.

This result raises the question whether exponential tail bounds hold for other relevant
ways of empirically (i.e. without using population parameters) standardizing sums.
The answer turns out to be positive 
and this paper establishes tail bounds that are sub-Gaussian or even closer to normal for the following settings:
Standardization by empirically centering and Studentizing sums of normal observations in Section~\ref{tscores}, 
standardization with the
log likelihood ratio statistic for observations from an exponential family in Section~\ref{loglik},
and standardization via self-normalization for observations from a symmetric distribution 
with unknown center of symmetry in Section~\ref{selfnorm}. The latter standardization give rise to a novel
scan statistic for heteroscedastic data that is based on self-normalization, and its asymptotic power
properties are also analyzed in Section~\ref{selfnorm}. This analysis shows that the tail bounds
are tight in the sense that they allow optimal detection in a certain scan problem; it is known that this 
optimality hinges on having the correct sub-Gaussian tail bound. 

\section{Normal tail bounds for Studentized constrasts and empirically centered sums}  
\label{tscores}

In order to derive a tail bound for empirically centered and Studentized sums it is convenient to establish
a more general result about Studentized linear contrasts:

\begin{Theorem}  \label{tstatistic2}
Let $X_1,\ldots,X_n$ i.i.d. N($\mu,\sigma^2$) and ${\bf b} \in \R^n$ with $\sum_{i=1}^n b_i=0$,
$\sum_{i=1}^n b_i^2 =1$. Then
$$
V\ :=\ \frac{\sum_{i=1}^n b_iX_i}{\sqrt{\frac{1}{n-1} \sum_{i=1}^n \bigl( X_i -\overline{X}\bigr)^2}}
$$ is a pivot and satisfies a normal tail bound:
\begin{align*}
V &\stackrel{d}{=} \frac{\sum_{i=1}^{n-1} Z_i}{\sqrt{\sum_{i=1}^{n-1} Z_i^2}} \qquad
\mbox{ for $Z_i$ i.i.d. N$(0,1)$},\\
\frac{V^2}{n-1} &\sim {\rm Beta}\left(\frac{1}{2},\frac{n-2}{2}\right),\\
\Pr (V>t) &\leq \Pr (\mbox{N$(0,1) >t)$ \quad for } 
\begin{cases}
t\geq 2.5 \mbox{ and } n\geq 10, \ \mbox{ or} \\
t\geq 2.75 \mbox{ and } n\geq 6,
\end{cases}
\end{align*}
and the analogous bound holds for the left tail of $V$.
\end{Theorem}

In particular, Theorem~\ref{tstatistic2} shows that the $t$-statistic regains the normal tail bound
if the location parameter is estimated in a suitable way. This follows by setting ${\bf b}=\frac{\bf c}{\sqrt{\sum_i c_i^2}}$
with $c_i:=1-\frac{m}{n}$
if $i \leq m$ and $c_i:=-\frac{m}{n}$ otherwise, which implies $\sum_i c_i =0$ and $\sum_i c_i^2=m \left(
1-\frac{m}{n}\right)$:

\begin{Corollary}  \label{tstatistic}
Let $X_1,\ldots,X_n$ i.i.d. N($\mu,\sigma^2$) and $\overline{X}=\frac{1}{n}\sum_{i=1}^n X_i$.
Then for $1 \leq m < n$:
$$
V\ :=\ \frac{\frac{1}{\sqrt{m\left(1-\frac{m}{n}\right)}} \sum_{i=1}^m \bigl(X_i -\overline{X}\bigr)}{
\sqrt{\frac{1}{n-1} \sum_{i=1}^n \bigl(X_i -\overline{X}\bigr)^2}}
$$
satisfies
$$
\Pr (V >t)\ \leq \ \Pr (\mbox{N$(0,1) >t)$ \quad for }
\begin{cases}
t\geq 2.5 \mbox{ and } n\geq 10, \ \mbox{ or} \\
t\geq 2.75 \mbox{ and } n\geq 6.
\end{cases}
$$
\end{Corollary}

Studentization is a special case of self-normalization, see e.g. de la Pe\~{n}a et al.~(2009)
and Section~\ref{selfnorm}. Self-normalization has certain advantages over standardizing with the
population standard deviation because, roughly speaking, erratic fluctuations of the statistic
are mirrored and therefore compensated by the random self-normalizing (Studentizing) term in the denominator,
see Shao and Zhou~(2016,2017) for formal results. 
Corollary~\ref{tstatistic} shows that centering empirically rather
than with the expected value can likewise be advantageous.
\medskip

{\sl Remark:} The algebraic tails of  the $t$-distribution can be bounded by an exponential bound if
the argument is small relative to the degrees of freedom, and this exponential tail bound may be useful for
certain applications that do not require bounds far out in the tails. A referee pointed out the following more general
example: If the $X_i$ are symmetric about $\mu$, then identity (1.1) in de la Pe\~{n}a et al.~(2009) gives
for $T$ in (\ref{Tstat}):
$$
\Pr (T >x)\ =\ \Pr \left(\frac{\sum_{i=1}^m (X_i-\mu)}{\sqrt{\sum_{i=1}^m (X_i-\mu)^2}} \geq \frac{ \sqrt{m} x}{
\sqrt{m-1+x^2}}\right) \ \leq \ \exp \left(-\frac{mx^2}{2(m-1+x^2)}\right)
$$
where the inequality follows from (\ref{Rademacher}). Hence $T$ has a sub-Gaussian tail for $x=O(\sqrt{m})$.
However, even for this restricted range of arguments this sub-Gaussian bound does not have the desired scale factor 1.
For example, $x=\sqrt{m}$ yields the bound $\exp (- x^2/(2c))$ with $c=2-\frac{1}{m}$, so even for large $m$ one 
obtains $c \approx 2$. The scale factor plays a key role in the theory and applications of sub-Gaussian tail bounds.

\section{Sub-Gaussian tail bounds for the log likelihood ratio statistic}
\label{loglik}

Let $X_1,\ldots,X_n$ be independent observations from a regular one-dimensional natural exponential
family $\{f_{\theta}, \theta \in \Theta\}$, i.e. $f_{\theta}$ has a density with repect
to some $\sigma$-finite measure $\nu$ which is of the form $f_{\theta}(x)= \exp (\theta x
-A(\theta))\,h(x)$ and the natural parameter space $\Theta =\{\theta \in \R:\,
\int \exp ( \theta x) h(x) \nu (dx) < \infty \}$ is open.

In order to derive good finite sample tail bounds in this setting, it turns out that it is useful
to standardize with the log likelihood ratio statistic rather than by centering and scaling.
In more detail, let $1\leq m <n$ and $\th_0 \in \Theta$. Then the generalized log likelihood ratio
statistic based on the observations $X_1,\ldots,X_m$ is
\begin{align}
\logLRm (\theta_0) &= \log \frac{\sup_{\theta \in \Theta} \prod_{i=1}^m f_{\theta} (X_i)}{
\prod_{i=1}^m f_{\theta_0} (X_i)} \nn \\
&= \sup_{\theta \in \Theta} \left( (\theta-\th_0) \sum_{i=1}^m X_i -m
\Bl(A(\theta)-A(\th_0)\Br)\right) \label{logLRI}
\end{align}

The MLE $\hat{\th}_m$ is defined as the argmax of (\ref{logLRI}) if the argmax exists. Note that
$\logLRm(\th_0)$ is always well defined whether $\hat{\th}_m$ exists or not.

$\logLRm(\th_0)$ represents a standardization of the sum $\sum_{i=1}^m X_i$ since by Wilk's theorem
$2\, \logLRm (\theta_0)$ is asymptotically pivotal $\chi_1^2$ if the population parameter is $\theta_0$.
The idea pursued in this section is that $\sqrt{2\, \logLRm (\theta_0)}$ is therefore approximately standard normal, 
and hence
it might be possible to establish a {\sl finite sample} sub-Gaussian tail bound. In the binomial case 
such a tail bound was indeed established by Rivera and Walther~(2013), see also 
Harremo\"{e}s~(2016) for bounds when $m=1$.
This section first extends the binomial bound
to the exponential family case and then addresses the case of empirical standardization where the typically unknown
$\th_0$ is replaced by the MLE.

It should be pointed out that while the square root of the log likelihood ratio does not commonly appear 
in the current literature, it has a history as a statistic for inference in exponential families.
Barndorff-Nielsen~(1986) calls ${\rm sgn}(\hat{\th}_m -\th_0) \sqrt{2\, \logLRm (\theta_0)}$, as well as its
empirically standardized counterpart below, the {\sl signed likelihood ratio statistic}. 
Rivera and Walther~(2013), Frick et al.~(2014) and K\"{o}nig et al.~(2020) use this statistic for detection
problems. An important advantage of working with this standardization 
is that it allows to make full use of the power of the Chernoff bound, as can be seen from the proof of
Theorem~\ref{logLRThm}(a). The resulting tail bound is therefore tighter than those obtained from various
relaxations of the Chernoff bound such as the Hoeffding or Bennett bounds.

Usually $\th_0$ is not known. Then an empirical standardization is obtained with the MLE $\hat{\th}_n$
substituted into the log likelihood ratio statistic for all the observations $X_1,\ldots,X_n$:

\begin{align}
\logLRmn &= \log \frac{\Bl(\sup_{\theta \in \Theta} \prod_{i=1}^m f_{\theta} (X_i)\Br)
\left(\sup_{\theta \in \Theta} \prod_{i=m+1}^n f_{\theta} (X_i)\right)}{
\sup_{\theta \in \Theta} \prod_{i=1}^n f_{\theta} (X_i)} \label{logLRIIc} \\
&= \sup_{\theta \in \Theta} \left( \theta \sum_{i=1}^m X_i -m A(\theta)\right) +
\sup_{\theta \in \Theta} \left( \theta \sum_{i=m+1}^n X_i -(n-m) A(\theta)\right)-
\sup_{\theta \in \Theta} \left( \theta \sum_{i=1}^n X_i -n A(\theta)\right). \nn
\end{align}

As an aside, this statistic can be interpreted as the generalized log likelihood ratio test statistic 
for testing a common $\th$
against different $\th$ for $X_1,\ldots,X_m$ and $X_{m+1},\ldots,X_n$. The standardization $V$ in
Corollary~\ref{tstatistic} has the same interpretation. In fact, if $f_{\th}$ is 
N$(\th,\sigma)$ with unknown mean $\th$ and known $\sigma$,
then one computes that $\sqrt{2\,\logLRmn}$ equals $V$ with the sample variance replaced
by $\sigma^2$ in the definition of $V$.

As another example, if the $X_i$ are Bernoulli with unknown parameter $p \in (0,1)$, then
the natural parameter for the exponential family is $\th=\log \frac{p}{1-p}$.
One computes that $\logLRmn$ equals
$$
m \left( \overline{X}_m \log \frac{\overline{X}_m}{\overline{X}} +(1-\overline{X}_m)
\log \frac{1-\overline{X}_m}{1-\overline{X}}\right) +(n-m) \left( \overline{X}_{m^c} \log
\frac{\overline{X}_{m^c}}{\overline{X}} +(1-\overline{X}_{m^c}) \log
\frac{1-\overline{X}_{m^c}}{1-\overline{X}}\right)
$$
where $\overline{X}_m:=\frac{1}{m}\sum_{i=1}^m X_i$, $\overline{X}_{m^c}:=\frac{1}{n-m}\sum_{i=m+1}^n X_i$
 and $\overline{X}:=\frac{1}{n}\sum_{i=1}^n X_i$.
This statistic was proposed as a scan statistic by Kulldorff~(1997) and, despite its lengthy form,
has been widely adopted for scanning problems in computer science and statistics, see e.g. 
Neill and Moore~(2004a,2004b) and Walther~(2010).

\begin{Theorem}  \label{logLRThm}
Let $X_1,\ldots,X_n$ be i.i.d. $f_{\th_0} \in \{f_{\theta}, \theta \in \Theta\}$, a regular 
one-dimensional natural exponential family, and let $1\leq m <n$.
Then for $x>0$:
\begin{itemize}
\item[(a)] $$
\Pr_{\th_0} \left( \sqrt{2\, \logLRm (\th_0)} > x\right)\ \leq \ 2 \exp\left( -\frac{x^2}{2}\right)
$$
\item[(b)] $$
\Pr_{\th_0} \left( \sqrt{2\, \logLRmn } > x\right)\ \leq\ 
\begin{cases}
%\Pr_{\th_0} \left( \sqrt{2\, {\rm logLR}_{I,I^c}(\theta_0) } > x\right)\ \leq\ 
(4+2x^2) \exp\left( -\frac{x^2}{2}\right) \\
(4+2e) \exp\left( -\frac{x^2}{2}\right) \qquad \mbox{ if } x\leq \left(nC\right)^{1/6}
\end{cases}
$$
for a certain constant $C$.
\end{itemize}
\end{Theorem}

The bounds can be divided by 2 if one considers the signed square-root for one-sided inference.
The proof of (a) proceeds by inverting the Cram\'{e}r-Chernoff tail bound as in Rivera and Walther~(2013),
where this technique is employed for the binomial case. The bounds in (b) do not quite match the
bound in (a) and the author has not been able to establish the simple
$2 \exp(-x^2/2)$ bound for (b). Simulations suggest that in fact an even better bound holds which
is closer to the standard normal bound, i.e. a bound that gains the factor $1/x$ on the sub-Gaussian bound as in
(\ref{bound318}). Establishing such a bound is a relevant open problem given its importance 
for scan statistics, see Walther and Perry~(2019) and the references therein.

%%%%%%%%%%%%%%%%%%%%%%%%%%%%%%%%%%%%%%%%%%%%%%%%%%%%%%%%%%%%%%%%%%%%%%%%%%%%%%%%%%%
%Notes:  The first inequality in (b) is not
%difficult to derive once (a) is established, but it does not quite give the desired bound $2 \exp(-x^2/2)$.
%However, this inequality is derived by replacing $\hat{\th}_{[n]}$ by $\th_0$ in the denominator of
%(\ref{logLRIIc}), which would seem to incur considerable slack.
%It is conjectured that (b) also allows
%a bound that is not worse than $2 \exp(-x^2/2)$. The author is able to prove such a sub-Gaussian bound
% only for a range of $x$ values, as given in the second inequality of (b).
% Walther (2010) establishes a subGaussian bound as in (b) for the permutation distribution in the Bernoulli setting. 
%%%%%%%%%%%%%%%%%%%%%%%%%%%%%%%%%%%%%%%%%%%%%%%%%%%%%%%%%%%%%%%%%%%%%%%%%%%%%%%%%%%%%%%%%%%%%%%%%%%%%%%

\section{Tail bounds for self-normalized and empirically centered sums of symmetric random variables}
\label{selfnorm}

The goal of this section is to extend the results for i.i.d. normal observations in 
Section~\ref{tscores} to a setting that allows heteroscedastic observations with not necessarily equal
expected values. Clearly, some additional assumption is necessary. The methodology proposed below
allows to treat the case of independent
(not necessarily identically distributed) observations having symmetric distributions
 with unknown and possibly different centers of symmetry.

It is informative to recapitulate the short and well known argument for establishing
a sub-Gaussian tail bound via self-normalization in the case where the center of symmetry is known 
to be zero, see e.g. de la Pe\~{n}a et al.~(2009):
If $X_1,\ldots,X_m$ are independent and symmetric about 0, then introduce i.i.d.
Rademacher random variables $R_1,\ldots,R_m$, $\Pr(R_1=1)=\Pr(R_1=-1)=\frac{1}{2}$, which are
independent of the $X_i$. Then $X_i \stackrel{d}{=}R_i X_i$ and hence for $t>0$:
\be  \label{Rademacher}
\Pr \left(\frac{\sum_{i=1}^m X_i}{\sqrt{\sum_{i=1}^m X_i^2}} >t \right) =
\Pr \left(\frac{\sum_{i=1}^m R_i X_i}{\sqrt{\sum_{i=1}^m X_i^2}} >t \right) =
\Ex \Pr \left(\frac{\sum_{i=1}^m R_i X_i}{\sqrt{\sum_{i=1}^m X_i^2}} >t \Big|X_1,\ldots,X_m \right)\leq 
\exp \left( -\frac{t^2}{2} \right)
\ee
by Hoeffding's inequality. Hence the sub-Gaussian tail bound is inherited from the Rademacher sum.
Sub-Gaussianity for self-normalized sums has been investigated in a number of papers. In the i.i.d. case,
Gin\'{e} et al.~(1997)
show that if the self-normalized sums are stochastically bounded (which always holds if the law of $X_i$
is symmetric), then they are uniformly sub-Gaussian for some scale parameter.
Also for the i.i.d. case, Shao~(1999) established asymptotic Cram\'{e}r-type large deviation results
under the assumption of a finite third moment. For independent but not necessarily identically distributed
$X_i$, Jing et al.~(2003) establish a Cram\'{e}r-type large deviation result under certain finite moment assumptions 
when $\Ex X_i=0$.
In the case where the distributions of the $X_i$ are symmetric about 0,
Efron~(1969, pp.\,1285--1288) suggested that it should be possible to lower the
sub-Gaussian tail bound (\ref{Rademacher}) to the normal tail $\Pr ({\rm N(0,1)} > t)$
in the usual hypothesis testing range $t>1.65$, but Fig.~1 in Pinelis~(2007) shows that the normal tail
is too small by a factor of at least 1.2 for certain $t \in (2,3)$. However,
recent remarkable results by Pinelis~(2012) and Bentkus and Dzindzalieta~(2015) show that the 
sub-Gaussian tail bound (\ref{Rademacher}) for the Rademacher sum can be improved upon to a 
bound of the order $\frac{1}{t} \exp (-\frac{t^2}{2})$, namely to a multiple of $\Pr ({\rm N(0,1)} > t)$  
where the multiple is at most 3.18 and is even close to 1 for large $t$.
This tail bound will then translate to the sum $\sum_{i=1}^m X_i$ after self-normalization
 via the above argument. This makes the use of the self-normalization very attractive in this
setting, cf. the remarks in Section~\ref{tstatistic}. 

The first aim of this section is to extend these results to the case where the center of symmetry
is unknown and may vary between the $X_i$. At first glance, this would appear to be a hopeless
undertaking since the above Rademacher argument depends crucially on the symmetry about zero.
However, there are observations available outside the summation window $X_1,\ldots,X_m$ which can
be used for an empirical standardization.
The idea is to construct an empirical centering which eliminates the
unknown center of symmetry from the symmetrization argument, or which at least results in certain
bounds on the center of symmetry. The second step then is to show that these bounds
still allow for nearly normal tails.

For simplicity of exposition it is assumed in the following
that $n=mp$ for integers $m\geq 1$ and $p\geq 2$. If $m$ is much smaller
than $n$, as is typically the case for scan problems, then this can always be arranged by discarding
a small fraction of the observations if necessary. The proposed empirical centering is given by
a linear transformation $\widetilde{\bf X} ={\bf A X}$, where the matrix ${\bf A}$ satisfies the conditions in
Proposition~\ref{prop}. One example of such an empirical centering is 

\be \label{power1}
\wt{X}_i\ :=\ X_i -\frac{1}{p-1} \sum_{j= m+(i-1)(p-1)+1}^{m+i(p-1)} X_j, \qquad i=1,\ldots,m
\ee

Corresponding to the linear tranformation ${\bf A}$ write $\wt{\mu}_i := \sum_{j=1}^n a_{ij}\mu_j$, where
$\mu_j$ is the center of symmetry of $X_j$. Note that it is not assumed that the $X_j$ have a finite
expected value. In the following, the subscript
$I$ denotes averaging over the index set $I:=\{1,\ldots,m\}$, so $\mu_I := \frac{1}{m} \sum_{i =1}^m \mu_i$
and $\mu_{I^c} := \frac{1}{n-m} \sum_{i=m+1}^n \mu_i$.

\begin{Proposition}  \label{prop}
Let ${\bf A}$
be a $m \times n$ matrix that has $p$ non-zero entries in each row and one non-zero
entry in each column, and these entries are 1 in columns $\{i: i \leq m\}$ and $\frac{-1}{p-1}$
in columns $\{i: i >m\}$. \footnote{This uniquely determines $A$ up to permutations
of the columns $\{i: i\leq m\}$ and permutations of the columns $\{i: i >m\}$.}

Let $X_i,\ i=1,\ldots,n$, be independent and symmetric about $\mu_i$ (so the $X_i -\mu_i$
need not be identically distributed).
\begin{itemize}
\item[(a)] If $\widetilde{\bf X} ={\bf A X}$, then the self-normalized sum of the $\widetilde{X}_i$
satisfies
$$
\frac{\sum_{i=1}^{m} \wt{X}_i}{\sqrt{\sum_{i=1}^{m} \wt{X}_i^2}}\ =\
\frac{n}{n-m} \frac{\sum_{i =1}^m (X_i-\overline{X})}{\sqrt{{\bf X}^T {\bf A}^T {\bf A X}}}\ =:\ T_m
$$
\item[(b)] If $\mu_I \leq \mu_{I^c}$ and $\mu_i=\mu_I$ for all
$i \leq m$ and $\mu_i =\mu_{I^c}$ for all $i >m$, then
\be  \label{bound318}
\Pr \left( T_m \geq t \right) \ \leq \ \min \Bl( 3.18, g(t)\Br) \,\Pr \Bl({\rm N(0,1)} > t\Br)
\ee
for all $t>0$, where $g(t):=1+\frac{14.11 \phi(t)}{(9+t^2) (1-\Phi(t))} \ra 1$ as $t \ra \infty$.
\item[(c)] If $\mu_I \leq \mu_{I^c}$ and (\ref{A}) or (\ref{Aprime}) hold, then the tail bound
(\ref{bound318}) holds for $t \in \left(0,\sqrt{m}K\right)$ for some ${K=K(v)>0}$.

Condition (\ref{A}) requires that the $\wt{\mu}_i$ don't vary much:
\be  \label{A}
\sum_{i =1}^m (\wt{\mu}_i -\wt{\mu}_I)^2 \leq \ v \sum_{i =1}^m \wt{\mu}_i^2
\qquad \mbox{ for some $v \in [0,1)$}
\ee
Condition (\ref{Aprime}) requires that the $\{\mu_i, i\leq m\}$ don't vary much 
and likewise for $\{\mu_i, i> m\}$:
\be \label{Aprime}
\left.\begin{aligned}
\frac{1}{m} \sum_{i =1}^m (\mu_i -\mu_I)^2 \\
\frac{1}{n-m} \sum_{i=m+1}^n (\mu_i -\mu_{I^c})^2
\end{aligned}
\right\} \leq v (\mu_I -\mu_{I^c})^2 \qquad \mbox{ for some } v\geq 0.
\ee
\item[(d)] The analogous inequalities to (b) and (c) hold for the left tail of $T_m$
if $\mu_I \geq \mu_{I^c}$.
\end{itemize}
\end{Proposition}

The proof of Proposition~\ref{prop} shows that the transformed $\wt{X}_i$ is symmetric about
$\wt{\mu}_i$ which may not equal zero. Nevertheless, the self-normalized sum of the $\widetilde{X}_i$
satisfies the normal tail bound (\ref{bound318}) if the $\mu_i$ satisfy the conditions given in (b) or (c).
(b) is a standard assumption for testing against an elevated mean on $I$, see Yao~(1993).
Note that $T_m$ is similar to the statistic $V$ used in Corollary~\ref{tstatistic} for the homoscedastic
case. Indeed, the proof of Theorem~\ref{tstatistic2}
 shows that $V$ is the self-normalized sum of ${\bf BX}$ for a certain $(n-1)\times n$
matrix ${\bf B}$.

\subsection{Scanning heteroscedastic observations having symmetric log-concave distributions}

As the statistic $T_m$ appears to be new, it is incumbent to demonstrate its utility with
an analysis of its power. To this end this section considers the scan problem where one observes independent
$X_i$, $i=1,\ldots,n$, where each $X_i$ has a distribution that is symmetric about some $\mu_i$ and log-concave,
i.e. $X_i$ has a density of the form $f(x)=\exp \,\phi_i (x-\mu_i)$, where $\phi_i: {\bf R} \ra [-\infty,\infty)$
is a concave function that is symmetric about 0. Special cases of log-concave distributions are the class of
normal distributions, where $\phi_i$ is a quadratic, the class of Laplace distributions, where $\phi_i$ is
piecewise linear, and more generally all gamma distributions with shape parameter $\geq 1$, all Weibull 
distributions with exponent $\geq 1$ and all beta distributions with both parameters $\geq 1$.
Log-concave distributions represent
an attractive and useful nonparametric surrogate for the class of Gaussian distributions in a range of problems 
in inference and modeling, see e.g. the review papers of Walther~(2009), Saumard and Wellner~(2014)
and Samworth~(2018). 

The goal of the scan problem under consideration here is to detect an elevated mean
$\mu_I > \mu_{I^c}$ on some interval $I=(j,k]$. Both the starting point $j$ and the length $|I|=k-j$
are unknown, likewise the $\mu_i$ and the distributions of the $X_i$, i.e. the functions $\phi_i$, are unknown.
Thus this is the general setting of Proposition~\ref{prop} with the additional assumption of log-concavity.
The log-concavity assumption allows to establish a result about the asymptotic detection power of the statistic $T_m$
that is similar to the homoscedastic normal case.

$T_m$ tests for an elevated mean on the interval $I=(0,m]$. It is straightforward
to analyze a different interval $I=(j,k]$, e.g. by applying $T_{k-j}$ to the rearranged data vector
$(X_{j+1},\ldots,X_n,X_1,\ldots,X_j)$. Denote this statistic by $T_I$.
Analyzing all possible intervals $I \subset (0,n]$ gives rise
to a multiple testing problem that is addressed by combining the corresponding $T_I$ into a scan
statistic. Walther and Perry~(2019) analyze several ways for combining the $T_I$ such that
optimal inference is possible, e.g. the Bonferroni scan. The use of that scan requires the availability
of a tail bound for the null distribution of $T_I$, such as (\ref{bound318}). The Bonferroni scan and the normal tail bound 
(\ref{bound318}) give $T_I$ a critical value of the form $\sqrt{2 \log \frac{n}{|I|}} +\kappa_{n,I}(\alp)$ 
with $\kappa_{n,I}(\alp) =O(1)$, which follows as in the proof of Theorem~2 in Walther and Perry~(2019).
Thus (\ref{powerT}) in the following  theorem shows that the Bonferroni scan based on the $T_I$ has
asymptotic power 1 if the assumptions of the theorem are met. These assumptions are discussed following the 
statement of the theorem.

\begin{Theorem} \label{power}
Let the $X_i$, $i=1,\ldots,n$, be independent with a log-concave distribution that is symmetric about some
$\mu_i$. Set $\sigma_i^2:={\rm Var }\, X_i$, $I:=(0,m]$, let
${\bf A}$ be the linear transformation (\ref{power1}) and write $\wt{\sigma}_i^2 := {\rm Var }\, \wt{X}_i$.
Assume the $\mu_i$ satisfy (\ref{A}) or (\ref{Aprime}).

If $\mu_I -\mu_{I^c} \geq \sqrt{\frac{(2+\eps_n) \sigma_I^2 R_I \log \frac{n}{|I|}}{|I|}}$
with $\eps_n \sqrt{\log \frac{n}{|I|}} \ra \infty$, $|I| \geq (\log n)^2$ and $R_I :=\frac{ \sum_{i \in I}
\wt{\sigma}_i^2}{\sum_{i \in I} \sigma_i^2}$, and if
\be \label{B}
\frac{\sigma_j^2}{\sigma_I^2} \ \leq \ S \sqrt{\max_{i \in I} (j-i)} \qquad \mbox{for all } j\in \{1,\ldots,n\}
\mbox{ and some } S>0,
\ee
then
\be \label{Rst}
R_I \ \leq \ 1+2S \sqrt{\frac{|I|^2}{n}}
\ee
and
\be \label{powerT}
\Pr \left(T_I \ >\ \sqrt{2 \log \frac{n}{|I|}} + O(1) \right) \ \ra 1\qquad (n \ra \infty).
\ee
\nin This result likewise holds for intervals $I=(j,j+m]$, $0\leq j\leq n-m$, by applying the theorem 
to $(X_{j+1},\ldots,X_n$, $X_1,\ldots,X_j)$.
\end{Theorem}

In order to compare the power of this scan statistic to an optimal benchmark, we first
 consider the special case where $X_i \sim {\rm  N}(\mu_i,\sigma^2$). 
For this special case of homoscedastic normal observations
it is known that there is a precise condition under which detection is possible with asymptotic power 1:
$\mu_I-\mu_{I^c} \geq
\sqrt{\frac{(2+\eps_n) \sigma^2  \log \frac{n}{|I|}}{|I|}}$, provided that $\eps_n$ does not go to zero
too quickly:
$\eps_n \sqrt{\log \frac{n}{|I|}} \ra \infty$. One the other hand, dedection is impossible if `$\sqrt{2+\eps_n}$'
is replaced by `$\sqrt{2-\eps_n}$'. Hence $\sqrt{2}$ measures the difficulty of the detection problem,
and the theory of that problem shows that it affects this difficulty as an {\sl exponent}.
This explains the efforts in the literature to approach
$\sqrt{2}$ as fast as possible, and the rates
$\sqrt{2\pm \eps_n}$ given above appear to be the currently best known rates. Attaining the factor $\sqrt{2}$
hinges on having the correct scale factor in the sub-Gaussian null distribution of the test statistic.
References and summaries
of these results are given in Walther and Perry~(2019) and Walther~(2022).

Theorem~\ref{power} shows that in the practically important range
$|I| \leq \sqrt{\frac{n}{\log n}}$, the Bonferroni scan based on the $T_I$ does indeed have asymptotic 
power 1 if $\mu_I-\mu_{I^c}$ exceeds the above detection threshold for the homoscedastic normal case,
since (\ref{Rst}) gives $R_I=1+o(\eps_n)$ and $\sigma_I^2=\sigma^2$ by homoscedasticity.
It is notable that this Bonferroni scan, which is designed to deal with heteroscedastic symmetric observations,
allows optimal detection in the special case of homoscedastic normal data. In fact,
Theorem~\ref{power} shows that
it achieves the detection boundary for the homoscedastic normal case already provided only the $\sigma_i,\ i \in I$, are
equal and the $\sigma_i$ outside $I$ don't grow too quickly, as required in (\ref{B}). 

If the data are
heteroscedastic, then Theorem~\ref{power} requires that $\sigma^2$ needs to be replaced by 
$\sigma_I^2 R_I$ $=\frac{1}{|I|}\sum_{i \in I} \wt{\sigma}_i^2$ in the lower bound for $\mu_I -\mu_{I^c}$.
It is beyond the scope of this paper to analyze whether this condition is optimal.

There appears to be not much literature about the scanning problem with heteroscedastic observations,
presumably because it is difficult to derive appropriate methodology. For example, the recent work of Enikeeva~(2018) 
considers the heteroscedastic Gaussian detection problem where $\sigma$ is allowed to be different on $I$ and
$I^c$, but it is assumed that $\sigma$ is constant and known on both $I$ and on $I^c$.
The finite-sample tail bound (\ref{bound318}) holds without such a restriction and thus self-normalized statistics 
may prove to be quite useful for scanning problems.

The proof of Theorem~\ref{power}  uses the following moment inequality for log-concave distributions, 
which may be of independent interest:

\begin{Proposition}  \label{moments}
If $X$ has a log-concave distribution that is symmetric about 0, then for all real numbers $r,s >0$:
$$
\Ex |X|^s \ \leq \ \left( \Ex |X|^r \right)^{\frac{s}{r}} \Gamma (s+1) (r+1)^{\frac{s}{r}}
$$
\end{Proposition}

If $0<s<r$, then $\Ex |X|^s \leq ( \Ex |X|^r )^{\frac{s}{r}}$ by Jensen's inequality, without any assumptions on the law of $X$.
The proposition shows that if the distribution is log-concave and symmetric, then it is possible to bound higher 
absolute moments in terms of lower absolute moments.

\section{Proofs}
\label{proofs}

\subsection{Proof of Theorem~\ref{tstatistic2}}

Write ${\bf X}=(X_1,\ldots,X_n)^T$ and
let ${\bf A}$ be an orthogonal $n \times n$ matrix with first row
$(\frac{1}{\sqrt{n}},\ldots,\frac{1}{\sqrt{n}})$. Then ${\bf Y:=AX}$ is a vector of independent
normal random variables with variance $\sigma^2$ and $\Ex Y_1=\sqrt{n} \mu$, $\Ex Y_i=0$, $i=2,\ldots,n$.
Further $\sum_{i=1}^n (X_i -\overline{X})^2 = \sum_{i=1}^n X_i^2 -n \overline{X}^2 = \sum_{i=1}^n
Y_i^2 -Y_1^2 = \sum_{i=2}^n Y_i^2$. Note that this is the same transformation that is commonly used in
textbooks to derive the distribution of Student's $t$-statistic.
In the latter case one is interested in $\sqrt{n}\, \overline{X} =Y_1$, which is independent of
$\sum_{i=2}^n Y_i^2$. In contrast, the condition $\sum_{i=1}^n b_i=0$ ensures that $\sum_{i=1}^n b_iX_i$
is a function of $(Y_2,\ldots,Y_n)$ only:
\be \label{V}
V\ =\ \frac{\langle {\bf b},{\bf X} \rangle}{\sqrt{\frac{1}{n-1} \sum_{i=1}^n \bigl( X_i -\overline{X}\bigr)^2}}\
=\ \frac{\langle {\bf b},{\bf A}^T{\bf Y} \rangle}{\sqrt{\frac{1}{n-1} \sum_{i=2}^n Y_i^2}}\
=\ \frac{\langle {\bf c},{\bf Y} \rangle}{\sqrt{\frac{1}{n-1} \sum_{i=2}^n Y_i^2}},
\ee
where ${\bf c=Ab}$ has $c_1=\sum_{i=1}^n \frac{1}{\sqrt{n}} b_i =0$ and thus $\sum_{i=2}^n c_i^2
=\sum_{i=1}^n c_i^2 = \sum_{i=1}^n b_i^2 = 1$.

Set $U_i:= \frac{Y_i}{\sqrt{\sum_{i=2}^n Y_i^2}}$, $i=2,\ldots,n$. Then ${\bf U}=(U_2,\ldots,U_n)^T$
has the uniform distribution on the $(n-2)$-dimensional unit sphere in $\R^{n-1}$ since the $Y_i$
are i.i.d. N$(0,\sigma^2)$. Therefore $\sum_{i=2}^n w_i U_i$, the length of the projection of ${\bf U}$
onto a unit vector ${\bf w}=(w_2,\ldots,w_n)^T$, has the same distribution for every unit vector ${\bf w}$.

Setting ${\bf w}=\left(\frac{1}{\sqrt{n-1}},\ldots,\frac{1}{\sqrt{n-1}}\right)^T$ gives\footnote{Alternatively,
construct rows 2 to $n$ of the orthogonal matrix ${\bf A}$ such that 
${\bf Ab=c}=\left(0,\sqrt{\frac{1}{n-1}},\ldots,\sqrt{\frac{1}{n-1}}\right)^T$. 
Then (\ref{V}) gives $V=\left(\sum_{i=2}^n Y_i\right)/\sqrt{\sum_{i=2}^n Y_i^2}$ without assuming that
the $X_i$ are normal. This also shows that $V$ is a self-normalized sum.
However, the $Y_i$ may not be independent if the $X_i$ are not normal.}
$$
V\ =\ \sqrt{n-1}\, \sum_{i=2}^n c_i U_i\ \stackrel{d}{=}\ \sqrt{n-1}\, \sum_{i=2}^n w_i U_i\
=\ \frac{\sum_{i=2}^n Y_i}{\sqrt{\sum_{i=2}^n Y_i^2}}\ \stackrel{d}{=}\ 
\frac{\sum_{i=1}^{n-1} Z_i}{\sqrt{\sum_{i=1}^{n-1}Z_i^2}}
$$
where the $Z_i$ are i.i.d. N$(0,1)$.
Setting ${\bf w}=(1,0,\ldots,0)^T$ gives
$$
V\ \stackrel{d}{=}\ \sqrt{n-1} \,\sum_{i=2}^n w_i U_i\ 
=\ \sqrt{n-1}\,\frac{Y_2}{\sqrt{\sum_{i=2}^n Y_i^2}}\ \stackrel{d}{=}\
\sqrt{n-1}\,\frac{Z_1}{\sqrt{\sum_{i=1}^{n-1} Z_i^2}},
$$
so $\frac{V^2}{n-1} \sim {\rm Beta}\left( \frac{1}{2},\frac{n-2}{2}\right)$ follows from a well
known fact about the beta distribution.

It is also known that the uniform distribution on the sphere in $\R^m$, $m:=n-1$, gives $U_2$
the density $\frac{\Gamma (\frac{m}{2})}{\Gamma (\frac{1}{2}) \Gamma (\frac{m-1}{2})}
\left( 1-u^2 \right)^{\frac{m-3}{2}} 1(u \in (-1,1))$,
hence $V \stackrel{d}{=} \sqrt{n-1} \,U_2$ has density
$$
f_V(t)\ =\ \frac{1}{\sqrt{m}}\, \frac{\Gamma (\frac{m}{2})}{\Gamma (\frac{1}{2}) \Gamma (\frac{m-1}{2})}
\left( 1-\frac{t^2}{m} \right)^{\frac{m-3}{2}} 1(-\sqrt{m} \leq t \leq \sqrt{m}).
$$
The plan is to show that $f_V(t)$ is not larger than the standard normal density $\phi (t) =\frac{1}{\sqrt{
2 \pi}} \exp (-\frac{t^2}{2})$ for $t$ large enough. Clearly $f_V(t) \leq \phi (t)$ for $t>\sqrt{m}$.
For $t \in (0,\sqrt{m})$ one has $\Gamma (\frac{m}{2}) \leq \Gamma (\frac{m-1}{2}) \sqrt{\frac{m}{2}}$
for $m>2$ by Gautschi's inequality, and $\log (1+x) \leq x -\frac{x^2}{2}$ for $x \in (-1,0)$:
\begin{align}
f_V(t) &\leq \frac{1}{\sqrt{2\pi}} \exp \left( \frac{m-3}{2} \log \left( 1-\frac{t^2}{m}\right)\right)\nn \\
&\leq \frac{1}{\sqrt{2\pi}} \exp \left( \frac{m-3}{2} \left( -\frac{t^2}{m} -\frac{t^4}{2m^2}\right)\right)\nn \\
&= \phi (t) \exp \left( \frac{3}{2m} t^2 -\frac{m-3}{4m^2} t^4\right) \nn \\
&\leq \phi (t) \qquad \mbox{ for }\ t^2 \geq \frac{6m}{m-3} \label{t1}
\end{align}
The condition is satisfied if e.g. $t \geq 3 $ and $m=n-1\geq 9$. Less conservative bounds obtain by
employing higher order terms for bounding $\log (1+x)$. For example, $\log (1+x) \leq x -\frac{x^2}{2}
+\frac{x^3}{3}$ for $x=-\frac{t^2}{m} \in (-1,0)$ yields
$$
f_V(t)\ \leq \ \phi (t) \exp \left( \frac{3}{2m} t^2 -\frac{m-3}{4m^2} t^4 -\frac{m-3}{6m^3}t^6\right).
$$
Dividing the argument in the exponent by $\frac{m-3}{2m^2} t^2$ shows that the argument is non-positive if
$$
\frac{3m}{m-3} -\frac{1}{2}t^2 -\frac{1}{3m} t^4 \ \leq \ 0
$$
and this  inequality holds for $t^2 \geq g(m):=\frac{3}{4}m \left( \sqrt{1+\frac{16}{m-3}}-1 \right)$.
One checks numerically that
\be \label{t3}
\max_{m \in \{5,\ldots,8\}} g(m) \leq 2.75^2,\ \ \max_{m \in \{9,\ldots,75\}} g(m) \leq 2.5^2.
\ee
Therefore $f_V(t) \leq \phi(t)$ follows for $t \geq 2.5$ and $m >75$ from (\ref{t1}),
for $t \geq 2.5$ and $9 \leq m \leq 75$ from (\ref{t3}), and for $t \geq 2.75$ and
$m \geq 5$ from these results together with (\ref{t3}). The last claim of the theorem now
obtains with $n=m-1$. $\hfill \Box$

\subsection{Proof of Theorem~\ref{logLRThm}}

The proof of (a) proceeds by inverting the Cram\'{e}r-Chernoff tail bound,
 as in Rivera and Walther~(2013) for the binomial case. $X_1$ has moment generating function
$\Ex_{\th_0} \exp (tX_1) = \exp \left( A(\th_0 +t) -A(\th_0)\right)$ for $\th_0 +t \in \Th$.
Markov's inequality gives for $x> \Ex_{\th_0} X_1$:
\begin{align*}
\Pr_{\th_0} \left( \frac{1}{m} \sum_{i=1}^m X_i >x\right) &\leq \inf_{t\geq 0}
\frac{\Ex \exp (t \sum_{i=1}^m X_i)}{\exp (t m x)}\\
&\leq \exp \left\{ - \sup_{t\geq 0, t+\th_0 \in \Th} \Bl(t m x -m \bigl(A(\th_0 +t)-A(\th_0))\bigr)
\Br)\right\} \\
& = \exp \left\{ - \sup_{\th \in \Th:\, \th \geq \th_0} m \Bl( (\th -\th_0)x -\bigl(A(\th) -A(\th_0)\bigr) 
\Br)\right\} \\
&= \exp \left\{ - \logLRm (x,\th_0) \right\}
\end{align*}
where $\logLRm (x,\th_0) := \sup_{\th \in \Th} m \Bl( (\th -\th_0)x-\left(A(\th)-A(\th_0)\right)\Br)$. This
conclusion used the fact that the sup over $\{\th \in \Th:\, \th \geq \th_0\}$ equals the sup over $\{\th \in \Th\}$
since convexity of $A$ yields 
\be \label{Lst}
(\th -\th_0)x - \left(A(\th)-A(\th_0)\right) \ \leq \ (\th -\th_0)x -(\th -\th_0)A' (\th_0) 
\ee
and the RHS is negative if $\th< \th_0$ and $x> \Ex_{\th_0} X_1 =A' (\th_0)$. The following claim will be proved below:
\be \label{Lstst}
\text{ The function $x \mapsto \logLRm (x,\th_0)$ is continuous and strictly increasing on $[\Ex_{\th_0} X_1,
\infty)\cap \mathcal{M}^0$}
\ee
where $\mathcal{M}$ denotes the convex hull of the support of $f_{\th_0}$. 
Analogously one shows that for $x< \Ex_{\th_0} X_1$:
$$
\Pr_{\th_0} \left( \frac{1}{m} \sum_{i=1}^m X_i <x\right) \ \leq \ \exp \bigl\{ -\logLRm(x,\th_0)\bigr\}
$$
and $x \mapsto \logLRm (x,\th_0)$ is continuous and strictly decreasing  on $(-\infty,\Ex_{\th_0} X_1]
\cap \mathcal{M}^0$. Together with \linebreak
$\logLRm(\Ex_{\th_0} X_1,\th_0)=0$, which follows from (\ref{Lst}) and 
$\Ex_{\th_0} X_1=A'(\th_0) \in \mathcal{M}^0$, one obtains
$$
\Pr_{\th_0} \left( \logLRm \Bl( \frac{1}{m} \sum_{i =1}^m X_i, \th_0 \Br) > t\right)
\ \leq \ 2 \exp (-t)
$$
for $t>0$ and claim (a) follows. It remains to prove (\ref{Lstst}).
This follows from Lemma~6.7 in
Brown~(1986) or from a general result in convex analysis to the effect that the Legendre transform $\phi(x):=
\sup_{\th \in \Th} (\th x -A(\th))$ satisfies $\phi'(x) = \arg\max_{\th \in \Th} (\th x -A(\th))=:\th(x)$
if $x \in \mathcal{M}^0$ (in which case the MLE $\th(x)$ exists uniquely and is given by $\th(x)=A'^{-1}(x)$
by exponential family theory) and $\phi''(x) =1/A''(\th(x)) =1/\text{Var}_{\th(x)} X_1 >0$
since the exponential family is minimal. Hence $\logLRm(x,\th_0)$ is differentiable wrt $x \in \mathcal{M}^0$
and
$$
\frac{d}{dx} \logLRm(x,\th_0)\ =\ m\,\bigl(\th(x)-\th_0\bigr)
$$
It was shown above that if $x> \Ex_{\th_0} X_1$, then the maximizer $\th(x)$ satisfies $\th(x)\geq \th_0 $.
Now (\ref{Lstst}) follows from $\frac{d}{dx} \th(x) =\phi''(x)>0$ for $x \in \mathcal{M}^0$.
Part (a) of the theorem is proved.

\bigskip

As for part (b), by the definition (\ref{logLRIIc})

\be \label{5.1}
\logLRmn \ \leq \ {\rm logLR}_{m,n}(\theta_0)\ =\ {\rm logLR}_I(\th_0) +{\rm logLR}_{I^c}(\th_0) 
\ee
where $I:=\{1,\ldots,m\}$ and $I^c:=\{m+1,\ldots,n\}$ and for an index set $J$ write
\begin{align*}
{\rm logLR}_J (\theta_0) &=\ \log \frac{\sup_{\theta \in \Theta} \prod_{i \in J} f_{\theta} (X_i)}{
\prod_{i \in J} f_{\theta_0} (X_i)}\\
\overline{X}_J &= \frac{1}{\# J} \sum_{i \in J} X_i\ \mbox{ and } \ \overline{X} = \frac{1}{n}\sum_{i=1}^n X_i 
\end{align*}
So ${\rm logLR}_I (\theta_0)=\logLRm(\th_0)$.
The proof of (a) established for $x>0$:
\be \label{5.2}
\Pr_{\th_0} \Bl(\logLRI(\th_0)>x,\, \overline{X}_I \geq \Ex_{\th_0} X_1 \Br)\ \leq \ \exp(-x)
\ee
and the same bound holds with $\overline{X}_I < \Ex_{\th_0} X_1$ in place of $\overline{X}_I 
\geq \Ex_{\th_0} X_1$ or with $I^c$ in place of $I$. (\ref{5.2}) shows that $\logLRI(\th_0)\ 1(
\overline{X}_I \geq \Ex_{\th_0} X_1) \stackrel{d}{\leq} E$, where $E \sim {\rm Exp}(1)$.

Since $\{X_i, i\in I\}$ and $\{X_i, i \in I^c\}$ are independent and stochastic order is
preserved under convolution, one gets
\be \label{5.3}
\logLRI(\th_0)\ 1(\overline{X}_I \geq \Ex_{\th_0} X_1) + {\rm logLR}_{I^c}(\th_0)\ 1(\overline{X}_{I^c} 
\geq \Ex_{\th_0} X_1) \ \stackrel{d}{\leq} \ R
\ee
where $R$ has the Erlang distribution with density $t e^{-t} 1(t>0)$. As (\ref{5.3}) holds for all
possible combinations of '$\geq$' and '$<$' in the indicator functions, the union bound gives
$$
\Pr_{\th_0} \Bl(\logLRI(\th_0) + {\rm logLR}_{I^c}(\th_0) >x\Bigr) \ \leq \ 4\, \Pr (R>x) 
\ =\ 4 \int_x^{\infty} te^{-t} dt\ =\ 4(1+x) e^{-x}.
$$
Now the first inequality in (b) follows with (\ref{5.1}).
\medskip

As for the second inequality, set $\alp:=\frac{\#I}{n}=\frac{m}{n}$. Then for $x>0$: 
\begin{align}
\Pr_{\th_0} \Bl( \overline{X}_I &- \overline{X}_{I^c} >x\Br) \ =\
\Pr_{\th_0} \left( (1-\alp)\sum_{i \in I}X_i -\alp \sum_{i \not\in I}X_i >\alp (1-\alp)nx \right) \nn\\
&\leq \inf_{t\geq 0} \frac{\Ex_{\th_0} \exp \left( (1-\alp) t \sum_{i \in I}X_i -\alp t \sum_{i \not\in I}X_i
\right)}{\exp \left(\alp (1-\alp) tnx \right)} \nn\\
&= \exp \left\{ -\sup_{t \geq 0} n \Bl( \alp (1-\alp) tx -\alp \left[ A(\th_0 +(1-\alp)t) -A(\th_0)
\right] -(1-\alp) \left[A(\th_0 -\alp t) -A(\th_0)\right]\Br)\right\} \label{L1}
\end{align}
One way to proceed from here would be via a Taylor series approximation of $A$ in order to derive
an exponential tail bound for $\alp (1-\alp) \frac{(\overline{X}_I- \overline{X}_{I^c})^2}{2\sigma_0^2}$
and likewise approximate $\logLRmn$ by this quantity. But these approximations will create notable slack
in the tail bound, while the proof in (a) shows that tight bounds are possible by employing a statistic
that conforms to the Cram\'{e}r-Chernoff bound. To this end define for $x\geq 0$
$$
\wt{\rm logLR}_n(x)\ :=\ \sup_{t \geq 0} n \Bl( \alp (1-\alp) tx -\left[\alp A(\th_0 +(1-\alp)t) 
+(1-\alp) A(\th_0 -\alp t) -A(\th_0)\right]\Br)
$$
Then (\ref{L1}) gives
\be \label{L2}
\Pr_{\th_0} \Bl( \wt{\rm logLR}_n \left(\overline{X}_I - \overline{X}_{I^c}\right) \geq x,
\overline{X}_I - \overline{X}_{I^c} \geq 0 \Br)\ \leq \ \exp(-x)
\ee
since $\wt{\rm logLR}_n(\cdot)$ is strictly increasing with $\wt{\rm logLR}_n(0)=0$ by
Jensen's inequality.

 The goal now is to show that $\left|\wt{\rm logLR}_n \left(\overline{X}_I - 
\overline{X}_{I^c}\right) -\logLRmn \right|$ is small relative to $\wt{\rm logLR}_n \left(\overline{X}_I -
\overline{X}_{I^c}\right)$. This is done with the following Proposition~\ref{LP}, which gives a
general result about the MLE in natural exponential families, and with Lemma~\ref{LL}.
In order to motivate part (b) of the following proposition, recall that the exponential family
$\{f_{\th} (x), \th \in \Th\}$ can alternatively be parameterized by its mean value, and the mapping
$\th \mapsto \Ex_{\th} X =A'(\th)$ is a homeomorphism between $\Th$ and ${\mathcal M}^0$, the interior of the 
convex hull of the support of $f_{\th_0}$, see e.g. Brown~(1986). 
The MLE $\hat{\th}$ is given by the solution of $A'(\hat{\th})=
\overline{X}$ if it exists. It may fail to exist if $\overline{X}$ falls on the boundary of ${\mathcal M}$.
For example, if a binomial$(n,p)$ experiment results in $n$ successes, then $\overline{X}=1$, but in
the natural parametrization the supremum of the likelihood is approached as the natural parameter
$\th =\log \frac{p}{1-p} \ra \infty$, so the MLE $\hat{\th}$ does not exist. This issue usually
becomes negligible in an asymptotic analysis of the MLE, but it has to be accounted for in a finite
sample statement. 
\begin{Proposition} \label{LP}
Let $X_1,\ldots,X_n$ be i.i.d. from a regular one-dimensional natural exponential family
$\{f_{\th}, \th \in \Th\}$. For $\th_0 \in \Th$ write $\mu_0 =\Ex_{\th_0} X_1$, $\sigma_0^2 =
{\rm Var}_{\th_0} X_1$, and ${\rm logLR}_{n}(\th_0)$ is defined in (\ref{logLRI}).
\begin{itemize}
\item[(a)] If the MLE $\hat{\th}$ exists, then
\begin{align*}
n \Bl(\frac{\overline{X} -\mu_0}{\sigma_0} \Br)^2 &\leq 2\, {\rm logLR}_{n}(\th_0)\ \frac{M}{\sigma_0^2}\\
n (\hat{\theta} -\th_0)^2 \sigma_0^2 &\leq 2\, {\rm logLR}_{n}(\th_0)\ \frac{M\sigma_0^2}{m^2}\\
n (\overline{X} -\mu_0)(\hat{\th} -\th_0) &\leq 2\, {\rm logLR}_{n}(\th_0)\ \frac{M}{m}
\end{align*}
where $m =\min_{\th \mbox{ between } \th_0 \mbox{ and }\hat{\th}} A''(\th)$,
$M=\max_{\th \mbox{ between } \th_0 \mbox{ and }\hat{\th}} A''(\th)$.
\item[(b)] Let $\del >0$ such that $[\th_0-\del, \th_0 +\del] \in \Th$ and set
$d_{\del}:=\min_{\th=\th_0 \pm \del} \bigl( (\th -\th_0) A'(\th)-(A(\th)-A(\th_0)\bigr)$.
Then $d_{\del}>0$. If ${\rm logLR}_{n}(\th_0) \leq n d_{\del}$, then the MLE $\hat{\th}$ exists
and satisfies $|\hat{\th}-\th_0| \leq \del$.
\end{itemize}
\end{Proposition}

{\bf Proof of Proposition~\ref{LP}:} As for part (a), Taylor's theorem gives for $\th$ between
$\th_0$ and $\hat{\th}$:
$$
A(\th) -A(\th_0) \ \leq\ A'(\th_0) (\th -\th_0)+\frac{M}{2} (\th -\th_0)^2
$$
Therefore these $\th$ satisfy
$$
\frac{{\rm logLR}_{n}(\th_0)}{n}\ \geq\ (\th -\th_0) \overline{X} -\bigl(A(\th)-A(\th_0)\bigr)
\ \geq \ (\th -\th_0)(\overline{X} -\mu_0) -\frac{M}{2} (\th -\th_0)^2
$$
as $A'(\th_0)=\mu_0$. Setting $\th:=\th_M:=\th_0 +\frac{\overline{X}-\mu_0}{M}$ one obtains
$$
\frac{{\rm logLR}_{n}(\th_0)}{n}\ \geq\ \frac{(\overline{X}-\mu_0)^2}{2M}
$$
provided it can be shown that 
\be \label{BW}
\th_M \mbox{ is between $\th_0$ and $\hat{\th}$.}
\ee
To this end, define the functions
\begin{align*}
L(\th) &:= (\th -\th_0) \overline{X} -\bigl( A(\th) -A(\th_0)\bigr)\\
g(\th) &:= (\th -\th_0)\left(\overline{X}-\mu_0 \right) -\frac{M}{2} (\th -\th_0)^2
\end{align*}
Then $L'(\th_0)=g'(\th_0) =\overline{X}-\mu_0$ and $L''(\th) \geq g''(\th)$ for $\th$ between
$\th_0$ and $\hat{\th}$. Hence one obtains $L'(\th) \geq g'(\th)$ for $\th \in [\th_0,\hat{\th}]$
if $\hat{\th} \geq \th_0$, and $L'(\th) \leq g'(\th)$ for $\th \in [\hat{\th},\th_0]$
if $\hat{\th} < \th_0$. 

Now consider the case $\hat{\th} \geq \th_0$. Since $A'(\hat{\th})=\overline{X}$ gives
$L'(\hat{\th})=0$, one gets $g'(\hat{\th}) \leq 0$. Since $\th_M$ is the maximizer of the quadratic
function $g(\th)$, $g'(\hat{\th}) \leq 0$ implies $\th_M \leq \hat{\th}$. 
\be  \label{2nd}
\overline{X}-\mu_0\ =\ A'(\hat{\th}) -A'(\th_0)\ =\ A''(\xi)(\hat{\th}-\th_0)\qquad \mbox{ for some
$\xi$ between $\th_0$ and $\hat{\th}$}
\ee
implies that $\overline{X}-\mu_0$ and $\hat{\th} -\th_0$ have the same sign, as $A''>0$.
So $\overline{X}-\mu_0 \geq 0$, but then $g(\th_0)=0$, $g'(\th_0)=\overline{X}-\mu_0 \geq 0$
and $g(\th_M) =\frac{(\overline{X}-\mu_0)^2}{2M}\geq 0$ implies $\th_0 \leq \th_M$. (If 
$\overline{X}-\mu_0=0$, then the quadratic $g$ has only one zero and $\th_M=\th_0$). This shows
(\ref{BW}) in the case $\hat{\th} \geq \th_0$, the case $\hat{\th} < \th_0$ is analogous.

The second inequality in (a) follows from (\ref{2nd}) which gives
$$
(\hat{\th}-\th_0)^2 \sigma_0^2\ \leq \ \frac{(\overline{X}-\mu_0)^2}{m^2} \sigma_0^2\ \leq\ 
2\,\frac{{\rm logLR}_{n}(\th_0)}{n} \frac{M\sigma_0^2}{m^2}
$$
as well as
$$
(\overline{X}-\mu_0) (\hat{\th}-\th_0) \ \leq \ \frac{(\overline{X}-\mu_0)^2}{m} \ \leq \
2\,\frac{{\rm logLR}_{n}(\th_0)}{n} \frac{M}{m}.
$$
\medskip

As for part (b), the function $h(\th):=(\th-\th_0)A'(\th) -\bigl(A(\th)-A(\th_0)\bigr)$
is stricly decreasing for $\th <\th_0$ and strictly increasing for $\th > \th_0$ since
$h'(\th)=(\th-\th_0) A''(\th)$. Further, $h(\th_0)=0$ and $\min_{\th=\th_0 \pm \delta} h(\th)
=d_{\delta}$. This shows that $d_{\delta}>0$ and
\be \label{b1}
\{\th \in \Th:\,h(\th) \leq d_{\delta}\}\ =:\ [\th_{low},\th_{up}]\ \subset\ [\th_0-\delta,
\th_0+\delta].
\ee
The motivation for defining $h$ is that for each $\th$, $h(\th)$ gives ${\rm logLR}_{n}(\th_0)/n$
when $\overline{X}=A'(\th)$, with $\th$ representing the argmax (i.e. the MLE). Indeed
\be  \label{b2}
h(\th)\ =\ \sup_{t\in \Th} \Bl[(t-\th_0) A'(\th) -\bigl(A(t)-A(\th_0)\bigr)\Br],\ \th \in \Th
\ee 
as is readily seen by differentiating wrt $t$. To make clear the dependence of 
${\rm logLR}_{n}(\th_0)$ on $\overline{X}$ we write similarly as before 
${\rm logLR}_{n}(\overline{X},\th_0):= {\rm logLR}_{n}(\th_0)$, i.e.
\be \label{b3}
\frac{1}{n} {\rm logLR}_{n}(x,\th_0)\ =\ \sup_{t \in \Th} \Bl[(t-\th_0)x-\bigl(A(t)-A(\th_0)\bigr)\Br],\ 
\ x \in {\mathcal M}.
\ee
This function is convex in $x$ since it is the Legendre transform of the convex function $A(t)$ plus a linear function.
Comparing (\ref{b2}) and (\ref{b3}) shows that
$$
\frac{1}{n} {\rm logLR}_{n}(x,\th_0)\ =\ h(\th)\qquad \mbox{ with } x=A'(\th),
$$
so this identity holds for $\th \in \Th$ and $x\in {\mathcal M}^0$, with $\th$ being the MLE
when the mean is $x$. Therefore
\be \label{b4}
\left\{ x\in {\mathcal M}^0:\ \frac{1}{n} {\rm logLR}_{n}(x,\th_0) \leq d_{\delta}\right\}\ =\ 
\bigl[A'(\th_{low}),A'(\th_{up})\bigr]
\ee
(recall that $A'$ is strictly increasing and continuous). But this implies that a boundary
point $x\in {\rm bd} {\mathcal M}$ cannot satisfy $\frac{1}{n} {\rm logLR}_{n}(x,\th_0) \leq
d_{\delta}$ because the function $x \mapsto {\rm logLR}_{n}(x,\th_0)$ is convex and hence
$M_{\delta}:=\ \bigl\{ x\in {\mathcal M}:\ \frac{1}{n} {\rm logLR}_{n}(x,\th_0)
\leq d_{\delta}\bigr\}$ must be an interval. Together with (\ref{b4}) this shows that
$M_{\delta} \in {\mathcal M}^0$ and so for every $x \in M_{\delta}$ the MLE exists and is
given by $(A')^{-1}(x) \in [\th_{low},\th_{up}] \subset [\th_0-\delta,\th_0+\delta]$.
$\hfill \Box$
\medskip

\begin{Lemma} \label{LL}
Let $T>0$. If the MLEs $\hat{\th}_I$ and $\hat{\th}_{I^c}$ exist, then on the event
$\bigl\{ \overline{X}_I -\overline{X}_{I^c} \geq 0,\, \sqrt{\alp} |\hat{\th}_I -\th_0| \leq
T,\, \sqrt{1-\alp} |\hat{\th}_{I^c} -\th_0| \leq T\bigr\}$:
$$
\logLRmn \ \leq \ \wt{\rm logLR}_n \left(\overline{X}_I - \overline{X}_{I^c}\right) +
\frac{5}{\sqrt{2\alp (1-\alp)}}\, nT^3\, \max_{\th:|\th-\th_0|\leq \frac{3T}{\sqrt{2 \alp (1-\alp)}}} |A'''(\th)|
$$
\end{Lemma}

{\bf Proof of Lemma~\ref{LL}:}
In the case where the MLEs $\hat{\th}_I$ and $\hat{\th}_{I^c}$ exist, set $\wt{\th}:=\alp
\hat{\th}_I +(1-\alp)\hat{\th}_{I^c}$. By definition (\ref{logLRIIc}):
\begin{align*}
\logLRmn &\leq {\rm logLR}_{m,n}(\wt{\th}) \\
&= \alp n \bigl( \hat{\th}_I \overline{X}_I -A(\hat{\th}_I)\bigr) +
(1-\alp) n \bigl( \hat{\th}_{I^c} \overline{X}_{I^c} -A(\hat{\th}_{I^c})\bigr) -n\bigl(\wt{\th}\overline{X}
-A(\wt{\th})\bigr) \\
&= \alp (1-\alp) n \bigl(\hat{\th}_I-\hat{\th}_{I^c}\bigr) \bigl(\overline{X}_I-\overline{X}_{I^c}\bigr)
-n \bigl[ \alp A(\wt{\th} +(1-\alp)t) -(1-\alp) A(\wt{\th} -\alp t) -A(\wt{\th}) \bigr]\\
& \qquad \mbox{ with $t:=\hat{\th}_I-\hat{\th}_{I^c}$ and using $\overline{X}=\alp \overline{X}_I
+(1-\alp) \overline{X}_{I^c}$} \\
&\leq \wt{\rm logLR}_n \bigl( \overline{X}_I-\overline{X}_{I^c}\bigr) + R\bigl(\hat{\th}_I-\hat{\th}_{I^c},
\wt{\th} \bigr)
\end{align*}
on $\bigl\{ \overline{X}_I-\overline{X}_{I^c}\geq 0 \bigr\}$, where
$$
R(t,\wt{\th})\ :=\ n\Bl[ \alp A(\th_0+(1-\alp)t) +(1-\alp) A(\th_0 -\alp t) -A(\th_0)\Br]
-n\Bl[ \alp A(\wt{\th}+(1-\alp)t) +(1-\alp) A(\wt{\th} -\alp t) -A(\wt{\th}) \Br].
$$
The last inequality uses the fact that $\hat{\th}_I-\hat{\th}_{I^c}$ and $\overline{X}_I-\overline{X}_{I^c}$
have the same sign since $\overline{X}_I-\overline{X}_{I^c} =A'(\hat{\th}_I)-A'(\hat{\th}_{I^c})$
and $A''>0$. 

Taylor's theorem gives for some $\xi,\tau$ between $0$ and $t$:
\begin{align*}
R(t,\wt{\th}) &= \frac{1}{2} \alp (1-\alp) n t^2 \Bl[ (1-\alp) A''(\th_0+(1-\alp) \xi) +\alp A''(\th_0 -\alp \xi)
-(1-\alp) A''(\wt{\th}+(1-\alp)\tau) -\alp A''(\wt{\th} -\alp \tau) \Br]\\
&\leq \frac{1}{2} \alp (1-\alp) n t^2 \Bl[ \max_{\th:|\th-\th_0|\leq |t|} A''(\th) -
\min_{\th:|\th-\wt{\th}|\leq |t|} A''(\th)\Br] \\
&\leq \frac{5}{\sqrt{2\alp (1-\alp)}}  n T^3 \max_{\th:|\th-\th_0|\leq \frac{3T}{\sqrt{2\alp (1-\alp)}}} |A'''(\th)|
\end{align*}
since $|\wt{\th}-\th_0| \leq \alp |\hat{\th}_I -\th_0| +(1-\alp) |\hat{\th}_{I^c}-\th_0| \leq
(\sqrt{\alp} +\sqrt{1-\alp})T \leq \sqrt{2}T$ and $|t|=|\hat{\th}_I-\hat{\th}_{I^c}| \leq
\left( \sqrt{\frac{1}{\alp}} +\sqrt{\frac{1}{1-\alp}}\right) T \leq \sqrt{\frac{2}{\alp (1-\alp)}} T$,
so $\left\{ \th: \max \left( |\th -\th_0|,|\th -\wt{\th}|\right) \leq |\hat{\th}_I -\hat{\th}_{I^c}|\right\}
\subset \left\{ \th: |\th -\th_0| \leq \frac{3T}{\sqrt{2\alp (1-\alp)}} \right\}$
and  $\max \left\{ |\th_1-\th_2|:\, |\th_1-\th_0| \leq |\hat{\th}_I -\hat{\th}_{I^c}|,\, 
|\th_2-\wt{\th}|\leq |\hat{\th}_I -\hat{\th}_{I^c}|\right\} \leq \frac{5T}{\sqrt{2\alp (1-\alp)}}$.
$\hfill \Box$
\bigskip

Now the proof of the theorem can be completed as follows: Let $\delta >0$ such that
$[\th_0-\delta, \th_0 +\delta] \subset \Th$. If $\logLRI(\th_0) \leq x$ for some 
$ x \in (0,\alp (1-\alp) n d_{\delta})$,
then part (b) of Proposition~\ref{LP} implies that the MLE $\hat{\th}_I$ exists and $|\hat{\th}_I-\th_0| \leq \delta$.
But then (a) of that Proposition implies that $\alp n (\hat{\th}_I-\th_0)^2 \leq 2x \frac{M}{m^2}$, where 
$M:= \max_{\th: |\th -\th_0| \leq \delta} A''(\th)$ and $m:= \min_{\th: |\th -\th_0| \leq \delta} A''(\th)$.
Likewise, ${\rm logLR}_{I^c}(\th_0) \leq x$ implies $(1-\alp) n (\hat{\th}_{I^c}-\th_0)^2 \leq 2x \frac{M}{m^2}$, hence
we can set $T:=\sqrt{\frac{2xM}{n m^2}}$ in Lemma~\ref{LL} to obtain on the event
$\{ \overline{X}_I-\overline{X}_{I^c} \geq 0,\, \logLRI \leq x,\, {\rm logLR}_{I^c}(\th_0) \leq x\}$:
$$
\logLRmn \ \leq \ \wt{\rm logLR}_n \bigl(\overline{X}_I-\overline{X}_{I^c} \bigr) +\sqrt{\frac{x^3}{n}}\,C
$$
where $C:= \frac{10}{m^3} \sqrt{\frac{M^3}{\alp (1-\alp)}} \max_{\th: |\th -\th_0| \leq \delta} |A'''(\th)|$.
So for $x \in \bigl( 0,n\, \min\bigl( \alp (1-\alp) d_{\delta}, C^{-2}\bigr)\bigr)$:
\begin{align*}
\Pr_{\th_0} &\left( \logLRmn >x,\, \overline{X}_I-\overline{X}_{I^c} \geq 0 \right) \\
&\leq \Pr_{\th_0} \left( \wt{\rm logLR}_n \bigl( \overline{X}_I-\overline{X}_{I^c} \bigr) > x\left( 1-\sqrt{\frac{x}{n}}
C\right),\, \overline{X}_I-\overline{X}_{I^c} \geq 0\right)\ +\ \Pr_{\th_0} \left(\logLRI (\th_0) >x\right)\\
&\qquad  +\Pr_{\th_0} \left({\rm logLR}_{I^c} (\th_0) >x\right)\\
&\leq \exp \left( -x\left(1-\sqrt{\frac{x}{n}}\,C\right)\right) + 2 \exp (-x)
 \qquad \mbox{ by (\ref{L2}) and part (a) of the theorem} \\
& \leq (2+e) \exp (-x) \qquad \mbox{ if } x \leq \left(n C^{-2} \right)^{1/3}.
\end{align*} 
The companion inequality with $\overline{X}_I-\overline{X}_{I^c} <0$ obtains analogously. The claim for
$\sqrt{ 2\, \logLRmn}$ follows for $\frac{1}{2}x^2 \leq (n C^{-2} )^{1/3}$, so one can use $8C^{-2}$
as the constant $C$ in the statement of the theorem.
$\hfill \Box$

\subsection{Proof of Proposition~\ref{prop}} 

The requirements for the matrix ${\bf A}$ imply that $\sum_{i=1}^m \wt{\mu}_i$ contains each $\mu_i$, $i\leq m$,
exactly once with coefficient 1,
and each $\mu_i$, $i >m$, exactly once with coefficient $\frac{-1}{p-1}$. Therefore
\be \label{matrixA}
\wt{\mu}_I\ =\ \frac{1}{m} \sum_{i=1}^m \wt{\mu}_i \ =\ \frac{1}{m} \sum_{i=1}^m \mu_i -
\frac{1}{m(p-1)} \sum_{i=m+1}^n \mu_i \ =\ \mu_I -\mu_{I^c}
\ee

\begin{Lemma}  \label{A1A2}
If (\ref{A}) or (\ref{Aprime}) hold, then
$$
\frac{m (\wt{\mu}_I)^2}{\sum_{i =1}^m \wt{\mu}_i^2}\ \geq \
\begin{cases}
1-v\ & \mbox{ if (\ref{A}) holds}\\
\frac{1}{4v+1} & \mbox{ if (\ref{Aprime}) holds.}
\end{cases}
$$
\end{Lemma}

\bigskip

\nin {\bf Proof of Lemma~\ref{A1A2}:} $\sum_{i =1}^m (\wt{\mu}_i -\wt{\mu}_I)^2 =\sum_{i =1}^m
\wt{\mu}_i^2 -m (\wt{\mu}_I)^2$ since $\wt{\mu}_I =\frac{1}{m} \sum_{i =1}^m \wt{\mu}_i$.
Hence (\ref{A}) bounds the RHS by $v \sum_{i =1}^m \wt{\mu}_i^2$, while (\ref{Aprime}) will
be shown to give the bound 
\be \label{lemmabound}
\sum_{i =1}^m (\wt{\mu}_i -\wt{\mu}_I)^2\ \leq \ 4 mv (\wt{\mu}_I)^2
\ee
so the claim follows in each case by collecting terms.

For simplicity of exposition (\ref{lemmabound}) will be proved for the linear transformation ${\bf A}$
given by (\ref{power1}). The proof goes through in the same way for a general matrix $A$ given
in Proposition~\ref{prop} by employing more cumbersome notation. Therefore
$\wt{\mu}_i  = \mu_i -\frac{1}{p-1} \sum_{j= m+(i-1)(p-1)+1}^{m+i(p-1)} \mu_j$ for $i=1,\ldots,m$.
Then it follows from (\ref{matrixA}) and Jensen's inequality that
\begin{align*}
\sum_{i =1}^m (\wt{\mu}_i -\wt{\mu}_I)^2 &= 4 \sum_{i =1}^m \left( \frac{\mu_i -\mu_I}{2}
-\sum_{j= m+(i-1)(p-1)+1}^{m+i(p-1)}\frac{\mu_j -\mu_{I^c}}{2(p-1)} \right)^2 \qquad \mbox{ since $\wt{\mu}_I=
\mu_I -\mu_{I^c}$} \\
&\leq 4 \sum_{i =1}^m \left( \frac{1}{2}(\mu_i -\mu_I)^2 +\frac{1}{2(p-1)}
  \sum_{j= m+(i-1)(p-1)+1}^{m+i(p-1)} (\mu_j -\mu_{I^c})^2 \right)\\
&= 2 \sum_{i =1}^m (\mu_i -\mu_I)^2 +\frac{2}{p-1} \sum_{j =m+1}^n (\mu_j -\mu_{I^c})^2 \\
&\leq 2 mv (\wt{\mu}_I)^2 + \frac{2}{p-1} (n-m) v (\wt{\mu}_I)^2 \qquad \mbox{ by (\ref{Aprime})}\\
&= 4mv (\wt{\mu}_I)^2 \qquad \mbox{ since $n-m=m(p-1)$.}
\end{align*}
$\hfill \Box$
\bigskip

As for proof of part (a) of the Proposition, by the construction of ${\bf A}$ the sum
$\sum_{i=1}^{m} \wt{X}_i$ contains each $X_i$, $i \leq m$, exactly once with coefficient 1,
and each $X_i$, $i >m$, exactly once with coefficient $\frac{-1}{p-1}$. Therefore
\begin{align*}
\sum_{i=1}^{m} \wt{X}_i &= \sum_{i =1}^m X_i -\frac{1}{p-1} \sum_{i=m+1}^n X_i \\
&=\frac{n}{n-m} \left( \sum_{i =1}^m \left(1-\frac{m}{n}\right) X_i -\frac{m}{n}
\sum_{i=m+1}^n X_i \right) \qquad \mbox{since } n=mp\\
&=\frac{n}{n-m} \sum_{i =1}^m (X_i -\overline{X}).
\end{align*}

As for (b) and (c), since each column of ${\bf A}$ has only one non-zero entry, it follows that if $i_1 \neq i_2$,
then $\wt{X}_{i_1}$ and $\wt{X}_{i_2}$ are functions of disjoint sets of $X_j$. Hence the
$\wt{X}_1,\ldots,\wt{X}_{m}$ are independent. 
$X_j -\mu_j \stackrel{d}{=} \mu_j -X_j$ and the independence of the $X_j$ yield
\begin{align*}
\Pr \left(\wt{X}_i -\wt{\mu}_i \leq t\right)\ &=\ \Pr \left(\sum_{j=1}^n a_{ij} (X_j -\mu_j) 
\leq t \right)
\ = \ \Pr \left(\sum_{j=1}^n a_{ij} (\mu_j -X_j) \leq t \right) \\
&= \Pr \left(\sum_{j=1}^n a_{ij} (X_j -\mu_j) \geq -t \right)
\ = \ \Pr \left(\wt{X}_i -\wt{\mu}_i \geq -t \right).
\end{align*}
Hence $\wt{X}_i$ is symmetric about $\wt{\mu}_i$.
Theorem~1.1 in Bentkus and Dzindzalieta~(2015) gives the bound 
$\frac{\overline{\Phi} (t)}{4 \overline{\Phi} (\sqrt{2})} \leq 3.18 \overline{\Phi}(t)$
for the self-normalized Rademacher sum and Theorem~1.1 in Pinelis~(2012) gives the bound
$\overline{\Phi} (t) +\frac{14.11 \phi (t)}{9+t^2}$. Hence the conditioning argument
(\ref{Rademacher}) yields
\be \label{318temp}
\Pr \left( \frac{\sum_{i=1}^{m} \left(\wt{X}_i-\wt{\mu}_i\right)}{\sqrt{\sum_{i=1}^{m}
\left(\wt{X}_i-\wt{\mu}_i\right)^2}} >t\right)\ \leq \ \min \Bl( 3.18, g(t)\Br) \,\Pr 
\Bl({\rm N(0,1)} > t\Br)
\ee
for all $t>0$, where $g(t):=1+\frac{14.11 \phi(t)}{(9+t^2) (1-\Phi(t))} \ra 1$ as $t \ra \infty$. 

Lemma~\ref{A1A2} gives
\be \label{Mdef}
\sqrt{ \sum_{i =1}^m \wt{\mu}_i^2}\ \leq \ M \sqrt{m}\, |\wt{\mu}_I|\qquad \mbox{ for some }M \geq 1.
\ee
Suppose $T_m=\frac{\sum_{i \leq m} \wt{X}_i }{\sqrt{\sum_{i \leq m} \wt{X}_i^2}} >t$ for some $t>0$.
Then $\sum_{i=1}^m (\wt{X}_i -\wt{\mu}_i) > 0$ since $\sum_{i=1}^m \wt{\mu}_i = m \wt{\mu}_I=
m(\mu_I - \mu_{I^c}) \leq 0$ by (\ref{matrixA}). Hence Minkowski's inequality gives
$$
\frac{\sum_{i =1}^m (\wt{X}_i -\wt{\mu}_i)}{\sqrt{ \sum_{i =1}^m (\wt{X}_i-\wt{\mu}_i)^2}}
\ \geq \ \frac{\sum_{i =1}^m \wt{X}_i -m \wt{\mu}_I}{\sqrt{ \sum_{i =1}^m \wt{X}_i^2}
+\sqrt{ \sum_{i =1}^m \wt{\mu}_i^2}}\ \geq\ \frac{t \sqrt{ \sum_{i =1}^m \wt{X}_i^2} +
m\, |\wt{\mu}_I|}{\sqrt{ \sum_{i =1}^m \wt{X}_i^2} +M\sqrt{m}|\wt{\mu}_I|} \ \geq\ t
$$
if $\sqrt{m} \geq Mt$. Hence for $ t\in \left(0,\sqrt{m}/M\right]$:
$$
\Pr \left( T_m > t\right) \ \leq \ \Pr \left(
\frac{\sum_{i =1}^m (\wt{X}_i -\wt{\mu}_i)}{\sqrt{ \sum_{i =1}^m (\wt{X}_i-\wt{\mu}_i)^2}} >t \right)
$$
and the last term satisfies (\ref{318temp}), proving (c).

If the $\mu_i$ are constant for $i\leq m$ and for $i>m$, then (\ref{Aprime}) holds with $v=0$,
so one can use $M=1$ in (\ref{Mdef}). Then $T_m$ satisfies the bound (\ref{bound318}) for all positive $t$ since
$\Pr (T_m >t) =0$ for $t>\sqrt{m}$ by Cauchy-Schwartz, proving (b).
(d) is analogous. $\Box$

\subsection{Proof of Theorem~\ref{power}:} 

On the event $E_n(I):= \left\{ \sum_{i \leq m} \wt{X}_i \geq 0\right\}$
Minkowski's inequality gives
$$
T_I \ \geq \ \frac{ \sum_{i \leq m} \wt{\mu}_i \ +\ \sum_{i \leq m} (\wt{X}_i -\wt{\mu}_i)}{
\sqrt{\sum_{i \leq m} \wt{\mu}_i^2} + \sqrt{ \sum_{i \leq m} (\wt{X}_i -\wt{\mu}_i)^2}}
$$
By Lemma~\ref{A1A2} there exists $M \geq 1$ such that
$$
\sqrt{\sum_{i \leq m} \wt{\mu}_i^2}\ \leq \ M \sqrt{m}\, |\wt{\mu}_I|\ =\ M\sqrt{m} (\mu_I -\mu_{I^c})
\qquad \mbox{ by (\ref{matrixA})}
$$
Set $Q_n(I) := \frac{\sum_{i \leq m} (\wt{X}_i -\wt{\mu}_i)}{
\sqrt{\sum_{i \leq m} \wt{\mu}_i^2} + \sqrt{ \sum_{i \leq m} (\wt{X}_i -\wt{\mu}_i)^2}}$
and $F_n(I):= \left\{ \sum_{i \leq m} (\wt{X}_i-\wt{\mu}_i)^2
\leq \sum_{i \leq m} \wt{\sigma}_i^2 (1+\frac{\eps_n}{4})\right\}$.
Then on the event $E_n(I) \cap F_n(I)$:
\begin{align*}
T_I &\geq \frac{m(\mu_I -\mu_{I^c})}{M\sqrt{m} (\mu_I -\mu_{I^c}) + \sqrt{\sum_{i \leq m}
\wt{\sigma}_i^2 \left(1+\frac{\eps_n}{4}\right)}} + Q_n(I) \qquad \mbox{ by (\ref{matrixA})}\\
&\geq  \frac{m\, \mu_{min} \bigl( \sum_{i \leq m} \wt{\sigma}_i^2 \bigr)^{-1/2}}{
M\sqrt{m}\, \mu_{min} \bigl( \sum_{i \leq m} \wt{\sigma}_i^2 \bigr)^{-1/2} +1+\frac{\eps_n}{8}} +Q_n(I)
\end{align*}

since $\mu_I-\mu_{I^c} \geq \mu_{min} := \sqrt{\frac{(2+\eps_n) \sigma_I^2 R_I \log \frac{n}{m}}{m}}$
and the function $x \mapsto \frac{ax}{bx+c}$ with $a,b,c>0$ is nondecreasing in $x>0$. Now $m\, \mu_{min}
\bigl( \sum_{i \leq m} \wt{\sigma}_i^2 \bigr)^{-1/2} = \sqrt{(2+\eps_n) \log \frac{n}{m}}$ since
$\sigma_I^2 R_I =m^{-1} \sum_{i \leq m} \sigma_i^2 R_I=m^{-1} \sum_{i \leq m} \wt{\sigma}_i^2$. 
Using $\frac{1}{1+y} \geq 1-y$ for $y>0$
one obtains on the event $E_n(I) \cap F_n(I)$:
\begin{align*}
T_I &\geq \sqrt{(2+\eps_n) \log \frac{n}{m}} \left(1-\frac{\eps_n}{8} -M \sqrt{\frac{
(2+\eps_n) \log \frac{n}{m}}{m}} \right) + Q_n(I) \\
&\geq \left( \sqrt{2} +\eps_n \left(\frac{1}{16} +o(1)\right)\right) \sqrt{\log \frac{n}{m}} +Q_n
\end{align*}
since $m \geq (\log n)^2$ and $(\log n)^{-1/2} =o(\eps_n)$.

Now 
$$
|Q_n(I)| \ \leq \ \frac{\left| \sum_{i \leq m} (\wt{X}_i -\Ex \wt{X}_i) \right|}{
\sqrt{ \sum_{i \leq m} (\wt{X}_i -\Ex \wt{X}_i)^2}}
$$
so both tails of $Q_n(I)$ satisfy the bound (\ref{bound318}) by (\ref{318temp}). 
Therefore $\Pr (T_I > \sqrt{ 2\log \frac{n}{m}}
+O(1)) \ra 1$ obtains (note that $\eps_n \sqrt{\log \frac{n}{m}}
\ra \infty$) once it is shown that $\Pr (E_n(I) \cap F_n(I) ) \ra 1$.

The proof of Proposition~\ref{prop} shows that the $\wt{X}_i$ are are independent and symmetric about
$\wt{\mu}_i$. Chebychev's inequality and $\sum_{i \leq m} \wt{\mu}_i=m
(\mu_I -\mu_{I^c})$ give
\begin{align*}
\Pr (E_n(I)^c ) &= \Pr \left( \sum_{i \leq m}(\wt{X}_i -\wt{\mu}_i) < -\sum_{i \leq m} \wt{\mu}_i \right)
 \ \leq \ \frac{\sum_{i \leq m} \wt{\sigma}_i^2}{m^2 (\mu_I -\mu_{I^c})^2} \\
&\leq \frac{\sum_{i \leq m} \wt{\sigma}_i^2}{m (2+\eps_n) \sigma_I^2 R_m \log \frac{n}{m}}
 \ =\ \frac{1}{(2+\eps_n) \log \frac{n}{m}} \ \ra 0\\
\Pr (F_n(I)^c ) &= \Pr \left(\sum_{i \leq m} \left( (\wt{X}_i -\wt{\mu}_i)^2 -\wt{\sigma}_i^2 \right)
 > \sum_{i \leq m} \wt{\sigma}_i^2 \frac{\eps_n}{4} \right) \\
&\leq 16 \frac{{\rm Var} \left( \sum_{i \leq m} (\wt{X}_i -\wt{\mu}_i)^2\right)}{
\left(\sum_{i \leq m} \wt{\sigma}_i^2 \right)^2 \eps_n^2} \ \leq \ 
16 \frac{\sum_{i \leq m} C (\wt{\sigma}_i^2)^2}{\left(\sum_{i \leq m} \wt{\sigma}_i^2 \right)^2 \eps_n^2}
\end{align*}
where $C:=\Gamma (5)\, 3^2 -1$ obtains by setting $r=2,s=4$ in Proposition~\ref{moments}. This uses the fact that 
 $\wt{X}_i -\wt{\mu}_i$ has a log-concave distribution since it is the sum of independent log-concave
random variables, see e.g. Saumard and Wellner~(2014). Now (\ref{B}) implies for $i \leq m$
\begin{align}
\wt{\sigma}_i^2\ =\ \sigma_i^2 + \frac{1}{(p-1)^2} \sum_{j = m+(i-1)(p-1)+1}^{m+i(p-1)} \sigma_j^2 
&\leq \sigma_i^2 + \frac{1}{(p-1)^2} (p-1) S \sqrt{n} \sigma_I^2 \nn \\
&\leq \sigma_i^2 +2S \sqrt{\frac{m}{p}} \sigma_I^2 \label{a} \\
&\leq 3S \sqrt{m} \sigma_I^2 \label{b} 
\end{align}
(\ref{b}) yields
$$
\Pr (F_n(I)^c)  \leq \ 16 C \frac{3 S \sqrt{m} \sigma_I^2}{\sum_{i \leq m} \wt{\sigma}_i^2
\eps_n^2} \ \leq \ 48 C \frac{S}{\sqrt{m} \eps_n^2} \ra 0
$$
since $m \geq (\log n)^2$ and $\eps_n \sqrt{\log n} \ra \infty$.

Finally, (\ref{Rst}) follows from (\ref{a}). \hfill $\Box$

\subsection{Proof of Proposition~\ref{moments}:}

The proof uses the following lemma repeatedly:

\begin{Lemma} \label{lemlast}
Let $w(x)$ be an integrable function on $(0,\infty)$ that does not change its sign from $+$ to $-$
as $x$ increases from 0 to $\infty$. Then $\int_0^{\infty} w(x)\,dx =0$ implies
$\int_t^{\infty} w(x)\,dx \geq 0$ for all $t>0$.
\end{Lemma}

The lemma obtains by observing that $\int_t^{\infty} w(x)\,dx <0$ for some $t>0$ implies $w(z)<0$
for some $z>t$ as well as $\int_0^t w(x)\,dx = \int_0^{\infty} w(x)\,dx - \int_t^{\infty} w(x)\,dx >0$,
which implies $w(s)>0$ for some $s \in (0,t)$, contradicting the assumption about the sign
changes of $w$.
\medskip

Since the density $f$ of $X$ is log-concave and symmetric about 0, it follows that $f$
is non-increasing on $[0,\infty )$ and that $f_0:=f(0) >0$. Hence
$$
f(x) \begin{cases}
  \leq u(x):=f_0\, {\bf 1} \left( |x| \leq \frac{1}{2f_0} \right) & \mbox{if } x \in \left(0,\frac{1}{2f_0} \right],\\
  \geq u(x) & \mbox{if } x > \frac{1}{2f_0}.
\end{cases}
$$
Set $w(x):=f(x)-u(x)$. Then $\int_0^{\infty} w(x)\,dx =\frac{1}{2}-\frac{1}{2}=0$ since both densities $f$ and $u$ are 
symmetric about 0. Hence Lemma~\ref{lemlast} gives
\be \label{moments1}
\int_t^{\infty} f(x) dx\ \geq \ \int_t^{\infty} u(x) dx\ \ \ \ \mbox{ for all } t>0.
\ee
Let $U \sim \mbox{Unif}\left(-\frac{1}{2f_0},\frac{1}{2f_0} \right)$. Then (\ref{moments1}) yields for $s>0$:
\begin{align}
\Ex |X|^s &= 2 \Ex |X|^s {\bf 1}(X>0)\ =\ 2 \int_0^{\infty} \Pr (X > t^{\frac{1}{s}}) dt \nn \\  
 &\geq 2 \int_0^{\infty} \Pr (U > t^{\frac{1}{s}}) dt\ =\  
 \Ex |U|^s\ =\ 2f_0 \int_0^{\frac{1}{2f_0}} u^s ds\ =\ \frac{(2f_0)^{-s}}{s+1} \label{moments2}
\end{align}
Let $V$ have density $v(x):=f_0 \exp \left(-2f_0 |x|\right)$. Since $\log f(x)$ is a concave function
and $\log v(x)$ is linear on $[0,\infty)$, the function $g(x):=\log v(x)-\log f(x)$ is convex
on $[0,\infty)$ and satsifies $g(0)=0$.  Hence $g(x)$ cannot change its sign from $+$ to $-$ as
$x$ increases from 0 to $\infty$, and this is therefore also true for $w(x):=v(x)-f(x)$. 
Again $\int_0^{\infty} w(x)\,dx =0$ holds since both densities $f$ and $v$ are symmetric about 0, so
Lemma~\ref{lemlast} gives
$$
\int_t^{\infty} f(x) dx\ \leq \ \int_t^{\infty} v(x) dx\ \ \ \ \mbox{ for all } t>0.
$$
Thus
\begin{align*}
\Ex |X|^s &= \ 2 \int_0^{\infty} \Pr (X > t^{\frac{1}{s}}) dt \nn \\
 &\leq 2 \int_0^{\infty} \Pr (V > t^{\frac{1}{s}}) dt\ =\
 \Ex |V|^s\ =\ 2f_0 \int_0^{\infty} v^s \exp (-2f_0 v) dv\ =\ \Gamma (s+1)\, (2f_0)^{-s} 
\end{align*}
Together with (\ref{moments2}) this shows that
$$
(s+1)^{-1} \ \leq \ (2 f_0)^s\, \Ex |X|^s \ \leq \ \Gamma (s+1) \ \ \ \ \mbox{ for all } s>0.
$$
Hence for any $r>0$:  $\Ex |X|^s \leq \Gamma (s+1) \left((2 f_0)^{-r}\right)^{s/r} \leq \Gamma (s+1) \left(
(r+1) \Ex |X|^r \right)^{s/r}$. \hfill $\Box$

\subsection*{References}

\begin{description}
\item Barndorff-Nielsen, O.E. (1986). Inference on full or partial parameters based on the
standardized signed log likelihood ratio. {\sl Biometrika} {\bf 73}, 307--322.
\item Bentkus, V.K. and Dzindzalieta, D. (2015) A tight Gaussian bound for weighted sums of
Rademacher random variables. {\sl Bernoulli} {\bf 21}, 1231--1237.
\item Boucheron, S., Lugosi, G., and Massart, P. (2013). {\sl Concentration inequalities: A nonasymptotic theory 
of independence}. Oxford University Press, Oxford, UK.
\item Brown, L.D. (1986). {\sl Fundamentals of Statistical Exponential Families}.
Institute of Mathematical Statistics, Hayward, CA.
\item[] de la Pe\~{n}a, V.H., Lai, T.L. and Shao, Q.M. (2009). {\sl Self-Normalized Processes:
Theory and Statistical Applications.} Springer, Berlin.
\item[] Enikeeva, F., Munk, A. and Werner, F. (2018). Bump detection in heterogeneous
Gaussian regression. {\sl Bernoulli} {\bf 24}, 1266--1306.
\item Efron, B. (1969). Student’s $t$-test under symmetry conditions. {\sl J. Amer. Statist. Assoc.}
{\bf 64}, 1278–-1302.
\item[] Frick, K., Munk, A. and Sieling, H. (2014). Multiscale change point inference.
{\sl J. R. Stat. Soc. Ser. B.} {\bf 76}, 495--580.
\item Harremo\"{e}s, P. (2016). Bounds on tail probabilities in exponential families. arXiv:1601.05179
\item[] Gin\'{e}, E., G\"{o}tze, F. and Mason, D. (1997). When is the Student t-statistic asymptotically
normal? {\sl Ann. Probab.} {\bf 25}, 1514-1531.
\item[] Jing, B. Y., Shao, Q. M. and Wang, Q. Y. (2003). Self-normalized Cram\'{e}r type large deviations
for independent random variables. {\sl Ann. Probab.} {\bf 31}, 2167–2215.
\item[] K\"{o}nig, C., Munk, A. and Werner, F. (2020). Multidimensional
multiscale scanning in exponential families: Limit theory and statistical consequences.
{\sl Ann. Statist.} {\bf 48}, 655-678.
\item[] Kulldorff, M. (1997). A spatial scan statistic. {\sl Comm. Statist. Theory Methods}
{\bf 26}, 1481–1496.
\item[] Neill, D. and Moore, A. (2004a). A fast multi-resolution method for detection of significant spatial 
disease clusters. {\sl Adv. Neural Inf. Process. Syst.} {\bf 10}, 651–658.
\item[] Neill, D. and Moore, A. (2004b). Rapid detection of significant spatial disease clusters. In {\sl Proc. 
Tenth ACM SIGKDD International Conference on Knowledge Discovery and Data Mining} 256– 265. ACM, New York.
\item Pinelis, I. (2007). Toward  the  best  constant  factor  for  the  Rademacher-Gaussian  tail comparison.
{\sl ESAIM Probab. Stat.} {\bf 11}, 412–-426.
\item Pinelis, I. (2012). An asymptotically Gaussian bound on the Rademacher tails.
{\sl Electron. J. Probab.} {\bf 17}, 1--22.
\item[] Rivera,  C. and Walther, G. (2013). Optimal detection of a jump in the intensity of a Poisson
process or in a density with likelihood ratio statistics. {\sl Scand. J. Stat.} {\bf 40}, 752-769.
\item Samworth, R.J. (2018). Recent progress in log-concave density estimation. {\sl Statist. Sci.} {\bf 33}, 493-509.
\item Saumard, A. and Wellner, J.A. (2014).  
Log-concavity and strong log-concavity: a review.
{\sl Statistics Surveys} {\bf 8}, 45-114.
\item Shao, Q.-M. (1999). Cram\'{e}r-type large deviation for Student's $t$ statistic. 
{\sl J. Theoret. Probab.} {\bf 12}, 387--398.
\item Shao, Q. and Zhou, W. (2016). Cram\'{e}r type moderate deviation theorems for self-normalized processes.
{\sl Bernoulli} {\bf 22}, 2029--2079.
\item Shao, Q. and Zhou, W. (2017). Self-normalization: Taming a wild population in a heavy-tailed world.
{\sl Appl.  Math.  J. Chinese Univ.} {\bf 32}, 253--269.
\item Shorack, G.R. and Wellner, J.A. (1986). {\sl Empirical Processes with Applications to Statistics}.
Wiley, New York.
\item van der Vaart, A. and Wellner, J. A. (1996). {\sl Weak Convergence and Empirical Processes. With
Applications to Statistics}. Springer, New York.
\item Vershynin, R. (2018). {\sl High-Dimensional Probability: An Introduction with Applications in Data
Science}. Cambridge University Press, Cambridge, UK.
\item Wainwright, M.J. (2019). {\sl High-Dimensional Statistics: A Non-Asymptotic Viewpoint}.
Cambridge University Press, Cambridge.
\item Walther, G. (2009). Inference and modeling with log-concave distributions.
{\sl Statist. Sci.} {\bf 24}, 319–327.
\item[] Walther, G. (2010). Optimal and fast detection of spatial clusters with scan statistics.
{\sl Ann. Statist.} {\bf 38}, 1010--1033.
\item Walther, G. and Perry, A. (2019). Calibrating the scan statistic: finite sample performance vs. asymptotics.
arXiv preprint arXiv:2008.06136.
\item Walther, G. (2022). Calbrating the scan statistic with size-dependent critical values: 
Heuristics, methodology and computation. In: Glaz, J, Koutras M.V. (eds)
{\sl Handbook of Scan Statistics}. Springer, New York, NY.
\item[] Yao, Q. (1993). Tests for change-points with epidemic alternatives.
{\sl Biometrika} {\bf 80}, 179--191.
\end{description}

\end{document}